\documentclass[12pt]{article}

\usepackage[english]{babel}
\usepackage{graphicx}
\usepackage{amsmath}
\usepackage{amsthm}
\usepackage{amsfonts}
\usepackage{amssymb}
\usepackage{mathabx}
\usepackage{epstopdf}
\usepackage{xcolor}
\usepackage{enumerate}
\newtheorem{defn}{Definition}[section]

\newtheorem{thm}[defn]{Theorem}

\newcommand{\be}{\begin{equation}}
\newcommand{\ee}{\end{equation}}
\newcommand{\bea}{\begin{eqnarray}}
\newcommand{\eea}{\end{eqnarray}}
\newcommand{\beas}{\begin{eqnarray*}}
\newcommand{\eeas}{\end{eqnarray*}}

\newcommand{\goto}{\rightarrow}

\newcommand{\bp}{\begin{proof}}
\newcommand{\ep}{\end{proof}}

\newcommand{\dstyle}{\displaystyle}

\begin{document}

\title{Extension to the Beraha-Kahane-Weiss \\Theorem with Applications}
\author{Jason Brown  \\ Department of Mathematics and Statistics, Dalhousie University\\ \\
\and 
Peter T. Otto  \\ Department of Mathematics,  Willamette University}

\maketitle



\begin{abstract}
\noindent The beautiful Beraha-Kahane-Weiss theorem has found many applications within graph theory, allowing for the determination of the limits of root of graph polynomials in settings as vast as chromatic polynomials, network reliability, and generating polynomials related to independence and domination. Here we extend the class of functions to which the BKW theorem can be applied, and provide some applications in combinatorics.
\end{abstract}


%
%

\section{Introduction} 

\medskip

There are many instances where {\em graph polynomials} arise. For example, the well known {\em chromatic polynomial} $\pi(G,x)$ (see, for example, \cite{dongbook}) counts the number of proper $x$-colourings of the vertex set of $G$, when $x$ is a nonnegative integer. The {\em (all-terminal) reliability polynomial} $\mbox{Rel}(G,p)$ is the function whose value  at $p \in [0,1]$ is the probability that the graph is connected, given that the vertices are always operational, but the edges are independently operational with probability $p$ (this model of reliability models the robustness of the network to random failures \cite{colbook}). There are as well many other related polynomials, including {\em two-terminal} \cite{colbook} and {\em strongly connected} \cite{browncox} reliability polynomials (the latter for directed graphs). Moreover, there are polynomials that are utilized in various graphical sequences. For example, the {\em independence polynomial} $i(G,x)$, is the generating function for the number of independent sets of each cardinality in the graph $G$; likewise, the {\em domination polynomial} serves a similar function for dominating sets of a graph.

In all cases, there has been a significant amount of interest in the zeros of such polynomials. The addiction of the Four Colour Problem (which can be succinctly stated that $4$ is \underline{never} the zero of a chromatic polynomial of a planar graph) led Tutte and others (in both the mathematics and theoretical physics arenas) to investigate the location of zeros of chromatic polynomials in the complex plane. An initial study of the zeros of reliability polynomials led to the conjecture in 1992 \cite{browncolbourn} that they all lay in the closed unit disk centered at $z = 1$, and it has only been proven false by the slimmest of margins \cite{roylesokal}. The zeros of each of chromatic, independence and domination polynomials have each been shown to be dense in the complex plane \cite{BHN,browntufts}  (while, of course, this has \underline{not} been shown for reliability polynomials).

Often graph polynomials for a family of graphs satisfy a fixed term recurrence 
\begin{eqnarray} 
P_{n+k}(x) & = & -\sum_{i=1}^{k}f_{i}(x)P_{n+k-j}(x), \label{BKW}
\end{eqnarray}
where the $f_{i}$'s are polynomials in $x$.
Such a recurrence can be solved using the usual method for linear recurrences to derive an explicit formula of the type
\begin{eqnarray} 
P_{n}(x) & = & \sum_{i=1}^{k} \alpha_i(x)(\lambda_{i}(x))^{n} \label{BKWexplicitformula}
\end{eqnarray}
where the $\lambda_i$'s are the zeros of the corresponding characteristic equation of the recursive relation \eqref{BKW}.
Beraha, Kahane and Weiss proved a beautiful result concerning the limits of the zeros of such polynomials. To be more precise, $z$ is a {\em limit of zeros} of the sequence of polynomials $P_{1},P_{2},\ldots$ if there is a sequence $z_{1},z_{2},\ldots$ of complex numbers such that $P_{n}(z_{n}) = 0$ and $\displaystyle{\lim_{n \rightarrow \infty}  z_{n} = z}$. Then the Beraha-Kahane-Weiss (BKW) Theorem is stated as follows:

\begin{thm}\label{Beraha}
Suppose that $P_{1},P_{2},\ldots$ satisfies (\ref{BKW}) and the following two {\em nondegeneracy conditions}:
\begin{itemize}
\item the polynomials $P_{n}$ do not satisfy a lower order recurrence than that in (\ref{BKW}), and 
\item there are no distinct $i$ and $j$ for which $\lambda_i = \omega \lambda_j$ for some $\omega$ of unit modulus. 
\end{itemize}
Then $z \in \mathbb{C}$ is a limit of zeros of the $P_{n}$ iff the $\lambda_{i}$ can be reordered such that one of the following holds:
\begin{enumerate}[i.]
\item $|\lambda_{k}(z)| > |\lambda_{i}(z)|$ for all $i \neq k$ and $\alpha_{k}(z) = 0$, or  
\item for some $l \geq 2$, $|\lambda_{1}(z)| = |\lambda_{2}(z)| = \cdots = |\lambda_{l}(z)| > |\lambda_{j}(z)|$ for all $j > l$.
\end{enumerate}
\end{thm}

\noindent The BKW theorem has been utilized to great effect for a number of graph polynomials (including Sokal's proof of the surprising result that chromatic roots are dense in the complex plane \cite{sokaldense}). Moreover, often, instead of a recurrence, one has a formula in the form of (\ref{BKWexplicitformula}), and one can apply the BKW Theorem as one can reverse the process and uncover the underlying recurrence (see \cite{BrownLimit}). In fact, Beraha, Kahane and Weiss' proof of their result is really of limits of zeros of sequences of polynomials whose form is given in (\ref{BKW}).

However, the statement of the BKW theorem is not quite accurate. When one solves a recurrence of the type in (\ref{BKW}), the zeros $\lambda_{i}$ of the characteristic polynomial of the recurrence may be repeated. Indeed, Beraha, Kahane and Weiss note this \cite{BKW}, stating that in such a case  (\ref{BKWexplicitformula}) ``is modified in the usual way, e.g., if $\lambda_{1}(z) = \lambda_{2}(z) \neq \lambda_{j}(z)$ for $j > 2$, the term $\alpha_{1}\lambda_{1}^{n}+ \alpha_{2}\lambda_{2}^{n}$ is replaced by $\alpha_{1}\lambda_{1}^{n}+ n\alpha_{2}\lambda_{2}^{n-1}$.'' And yet, in spite of the explicit mention in \cite{BKW}, a closer look reveals that repeated $\lambda_i$'s are not covered by the theorem and proof. There are a number of instances of formulas for graph polynomials that satisfy a recurrence with exactly such repeated zeros -- equivalently, the explicit formula has the $\alpha_{i}$ a polynomial function of both $z$ and $n$. In light, it is worthwhile to revisit the statement and proof of the BKW Theorem to extend it to such cases. And we shall apply the result to a variety of graph polynomials.

Furthermore, the impetus for our extension of the BKW Theorem began with a new graph polynomial, which we refer to as a {\it Steele polynomial} $S(G; t)$ of a graph $G$.  In \cite{Steele}, Steele derived an integral formula for the mean length of minimal spanning trees of $G$ with random edge lengths distributed uniformly over the unit interval $[0,1]$.  Specifically, if we let $L(G)$ denote the random total length of a minimal spanning tree of the graph $G$, then $\mathbb{E}[L(G)] = \int_0^1 S(G; t) \, dt$ where 
\be
\label{SteelePoly} 
S(G; t) = \frac{(1-t)}{t} \frac{T_x(G; 1/t, 1/(1-t))}{T(G; 1/t, 1/(1-t))}
\ee
and $T(G; x, y)$ is the well known Tutte polynomial of the graph $G$.  In \cite{NOS}, it was shown that $S(G; t)$ is a polynomial of degree at most equal to the number of edges of $G$ (additional properties of Steele polynomials, including information on the coefficients of the polynomials, appears in \cite{NOS}).

The Steele polynomial for the sequence of cycle graph $C_n$ has the form 
\be
\label{eqn:SteeleCn}
S(C_n; t) = t^n - n t + (n-1)
\ee
which satisfy the following recurrence relation
\[ P_n(t) = (t+2)P_{n-1}(t) - (2t+1)P_{n-2}(t) + tP_{n-3}(t). \]
The corresponding characteristic polynomial has a simple zero at $\lambda_1(t) = t$ and a repeated zero at $\lambda_2(t) = 1$ which (with appropriate initial conditions) yield the solution $S(C_n; t)$ above.  It is the repeated zero $\lambda_2(t) = 1$ that forces the Steel polynomials of cycle graphs beyond the scope of the original BKW Theorem (Theorem \ref{Beraha}).

\section{Extending the Baraha-Kahane-Weiss Theorem} 

In the next theorem, we state our extension to the BKW Theorem that includes the case where the zeros of the characteristic equation $\lambda_i$'s of the recursive relation \eqref{BKW} are not necessarily distinct but can be repeated roots of arbitrary order.  It is straightforward to see that this is equivalent to the condition that the coefficients $\alpha_i$'s in \eqref{BKWexplicitformula} can be functions of both $n$ and $x$.  This is the condition we state our extension below.

We state the result for two zeros $\lambda_1(x)$ and $\lambda_2(x)$, but the result can be extended to more terms, by a direct generalization of our proof of part (i) and as in the original BKW theorem for part (ii). 


\begin{thm}
\label{MainThm}
Let $\{ P_n(x) \}$ be a sequence of analytic functions of the form
\be
\label{eqn:P_n}
P_n(x)=\alpha_1(n; x)(\lambda_1(x))^n+\alpha_2(n; x)(\lambda_2(x))^n,
\ee
where $\lambda_i$ are analytic and nonzero such that $\lambda_1(x) \neq \omega\lambda_2(x)$ for any $\omega\in\mathbb{C}$ of unit modulus, and $\alpha_i(n; x)$ have the form
\be
\label{eqn:alphas}
\alpha_i(n; x) = n^{\alpha_i} p_{i, \alpha_i}(x) + n^{\alpha_i - 1} p_{i, \alpha_1-1}(x) + \cdots + n p_{i, 1}(x)+ p_{i, 0}(x).
\ee
where the coefficient functions $p_{i, j}$ are analytic and $p_{i, \alpha_i}$ are nonzero. 

Then $z\in\mathbb{C}$ is a limit of zeros of the family $\{P_n(x)\}$ if either of the following conditions hold.
\begin{enumerate}[i.]
\item $|\lambda_1(z)| > |\lambda_2(z)|$ and $p_{1, \alpha_1}(z)=0$, or $|\lambda_2(z)| > |\lambda_1(z)|$ and $p_{2, \alpha_2}(z)=0$.
\item $\left| \lambda_1(z) \right| = \left| \lambda_2(z) \right| > 0$ and at least one of $p_{1, \alpha_1}(z)$ and $p_{2, \alpha_2}(z)$ is nonzero.
\end{enumerate}
\end{thm}
\bp

The Beraha-Kahane-Weiss Theorem was proved using Rouch\'e's Theorem, and we shall make repeated use of this well known result  in the following form:
Suppose that two functions $f(z)$ and $g(z)$ are analytic inside and on a simple closed curve $C$.  If $|f(z)| > |g(z)|$ at each point $z$ on $C$, then $f(z)$ and $f(z) + g(z)$ have the same number of zeros inside $C$.

\medskip

\noindent {\bf Part (i):}  By symmetry, we assume that $|\lambda_1(z)| > |\lambda_2(z)|$ and $p_{1, \alpha_1}(z)=0$ and prove that $z$ is a limit of zeros. To do so, we base our proof on the approach taken for the corresponding case in \cite{BKW}.  We will show that for every $\epsilon > 0$ there is a sufficiently large $N$ such that for all $n \geq N$, $P_n(z_n) = 0$ for some $z_n \in D_\epsilon(z) = \{x : |x-z| < \epsilon \}$.  Let $C$ denote the boundary of $D_\epsilon(z)$.


Since the zeros of the nonzero analytic function $p_{1, \alpha_1}$ are isolated, by taking $\epsilon$ sufficiently small, we can assume that $|p_{1, \alpha_1}(x)| > 0$ for $x \in C$; moreover, by the condition on $\lambda_1(z)$ and $\lambda_2(z)$, we can assume as well that $|\lambda_1(z)/\lambda_2(z)| \leq \rho < 1$ on $C$ as well.
Define
\[ \overline{p}_{1, \alpha_1} (n; x) = \frac{\left(n^{\alpha_1-1}p_{1, \alpha_1 - 1}(x) + \cdots + n p_{i, 1}(x)+ p_{i, 0}(x)\right)}{n^{\alpha_1}}.\] 
Then 
\[ \left(p_{1, \alpha_1}(x) + \overline{p}_{1, \alpha_1} (n; x) \right) n^{\alpha_1} = \alpha_1(n; x). \]
Next define
\[ w_n(x) = - \left(\overline{p}_{1, \alpha_1} (n; x) + \frac{\alpha_2(n; x) \lambda_2(x)^n}{n^{\alpha_1} \lambda_1(x)^n} \right). \]
Since for $x$ on the boundary $C$, $|\overline{p}_{1, \alpha_1} (x)| \goto 0$ as $n \goto \infty$, $|\lambda_2(x)|/|\lambda_1(x)|  < 1$ and $\alpha_2(n; x)$ is polynomial in $n$, we conclude that $|w_n(x)|$ can be made sufficiently small on $C$ for $n$ large. Specifically, there exists $N > 0$ such that for all $n \geq N$, $|w_n(x)| < |p_{1, \alpha_1}(x)|$ for $x \in C$.  By Rouch\'e's Theorem, for $n \geq N$, $p_{1, \alpha_1}(x) - w_n(x)$ and $p_{1, \alpha_1}(x)$ have the same number of zeros in $D_\epsilon(z)$.  
Finally, since $p_{1, \alpha_1}(z)=0$, for all $n \geq N$, there exists at least one $z_n \in D_\epsilon(z)$ such that $p_{1, \alpha_1}(z_n) = w(z_n)$, which is equivalent to $P_n(z_n) = 0$.

%

\medskip

\noindent {\bf Part (ii):}  We base our proof again on that of \cite{BKW}.  From \cite{BKW}, we know that points satisfying $\left| \lambda_1 \right| = \left| \lambda_2 \right|$ are not isolated; that is, this equality holds for points arbitrarily close to $z$.  

Let $U$ be a disk about $z$ and define $\mu(x) = \lambda_1(x)/ \lambda_2(x)$.  The assumption that $\lambda_1(x) \neq \omega\lambda_2(x)$ for any $\omega\in\mathbb{C}$ of unit modulus implies that $\mu$ is not constant and the non-isolated nature of the points satisfying $\left| \lambda_1 \right| = \left| \lambda_2 \right|$ implies that we can assume that $\mu'(z) \neq 0$.  Thus, for $U$ sufficiently small, $\mu$ is invertible from $U$ onto a neighborhood $V$ of $\omega_0 = \mu(z)$.  By the assumption that $\left| \lambda_1(z) \right| = \left| \lambda_2(z) \right|$, we have $|\omega_0| = 1$. 

By symmetry, we can assume that $p_{1, \alpha_1}(z) \neq 0$, and hence for all sufficiently large $n$, $\alpha_1(n; z) \neq 0$ on $U$. If we denote $\nu$ to be the inverse of $\mu$, then
\be
\Pi_n(x) = \frac{P_n(v(x))}{\alpha_1(n; \nu(x)) \lambda_2(\nu(x))^n} = x^n + \frac{\alpha_2(n; \nu(x))}{\alpha_1(n; \nu(x))} 
\ee 
By taking a smaller (noncircular) neighborhood in $U$, we can assume that for some $r_0 > 1$ and $\theta_0 > 0$, the image of \[ V = N_{r_0, \theta_0} (\omega_0) = \{ re^{i\theta} \omega_0 : r_0^{-1} < r < r_0, -\theta_0 < \theta < \theta_0 \} \]
under $\nu$ is a subset of $U$.  We will use the observation that for $V$ small, the ratio $\frac{\alpha_2(n; v(x))}{\alpha_1(n; v(x))}$ will be close to $\frac{\alpha_2(n; z)}{\alpha_1(n; z)}$.

Let $C = C_1 + C_2 + C_3 + C_4$, where $|x| = r_0$ on $C_1$, $|x| = r_0^{-1}$ on $C_3$, $x = re^{i\theta_0} \omega_0$ on $C_2$, and $x = re^{-i\theta_0} \omega_0$ on $C_4$.  We traverse around $C$ in the counter-clockwise direction.  Since $\alpha_1$ and $\alpha_2$ are polynomial in $n$ and the modulus of $x$ on $C_1$ is greater than $1$, $x^n$ dominates $\frac{\alpha_2(n; v(x))}{\alpha_1(n; v(x))}$ for large $n$ on $C_1$ and the image of $C_1$ winds around $0$ on the order of $n\theta_0$ times.  On the other hand, along $C_3$ the modulus is less than $1$ and thus $\frac{\alpha_2(n; v(x))}{\alpha_1(n; v(x))}$ dominates and since $\frac{\alpha_2(n; v(x))}{\alpha_1(n; v(x))}$ is close to $\frac{\alpha_2(n; z)}{\alpha_1(n; z)}$, for $n$ large, the image of $C_3$ does not wind around $0$.

Lastly, along $C_2$ and $C_4$, the argument is fixed and $\frac{\alpha_2(n; v(x))}{\alpha_1(n; v(x))}$ is close to $\frac{\alpha_2(n; z)}{\alpha_1(n; z)}$, the images of these two curves also do not wind around $0$.  Therefore, we conclude that, for all $n$ sufficiently large, the image of $\Pi_n$ of the boundary $C$ of the set $V$ winds around $0$ a positive number of times.  This implies that for all $n$ sufficiently large, $\Pi_n$ has a zero in $V$ which implies that $P_n$ has a zero in $U$.  
\ep

\vspace{0.25in}
To derive a partial converse to the previous theorem,  suppose that $\displaystyle{\lim_{n \rightarrow \infty} z_n = z}$ with $P_n(z_n) = 0$, and without loss $|\lambda_1(z)| > |\lambda_2(z)|$.  Then we have 
\[ \alpha_1(n; z_n) \lambda_1(z_n)^n = - \alpha_2(n; z_n) \lambda_2(z_n)^n.\] 
which implies that 
\begin{eqnarray}
\left|\alpha_1(n; z_n)\right|^{1/n} \left| \lambda_1(z_n) \right| & = & \left| \alpha_2(n; z_n) \right|^{1/n} \left| \lambda_2(z_n) \right|. \label{converseBKWextension}
\end{eqnarray}
If $p_{1, \alpha_1}(z) \neq 0$ and $p_{2, \alpha_2}(z) \neq 0$ then close to $z$, both $\alpha_1(n; z_n)$ and $\alpha_2(n; z_n)$ are bounded away from $0$, and hence their $n$--th roots go to $1$. It follows by taking limits in (\ref{converseBKWextension}) that $|\lambda_{1}(z)| = |\lambda_2(z)|$, a contradiction.

%
%

\vspace{0.25in}

We remark that  while the form of the coefficient functions $\alpha_i$ in \eqref{eqn:alphas} as polynomials in $n$ cover the applications included in this work, they are sufficient but not necessary assumptions for the conclusion of Theorem \ref{MainThm}.  For example, for part (i), the assumption on the $\alpha_i$ could be generalized to 
\[ \left| \alpha_2(n; z) \right| = {\dstyle o \left( \left|\frac{\lambda_1(z)}{\lambda_2(z)} \right|^n \right)}. \]
For part (ii), the generalization extends to 
\[ \left| \frac{\alpha_2(n; v(x))}{\alpha_1(n; v(x))} \right| = {\dstyle o \left( \left|x \right|^n \right)}. \]

\medskip

\begin{figure}
\label{example2}
\begin{center}
 \includegraphics[scale=.8]{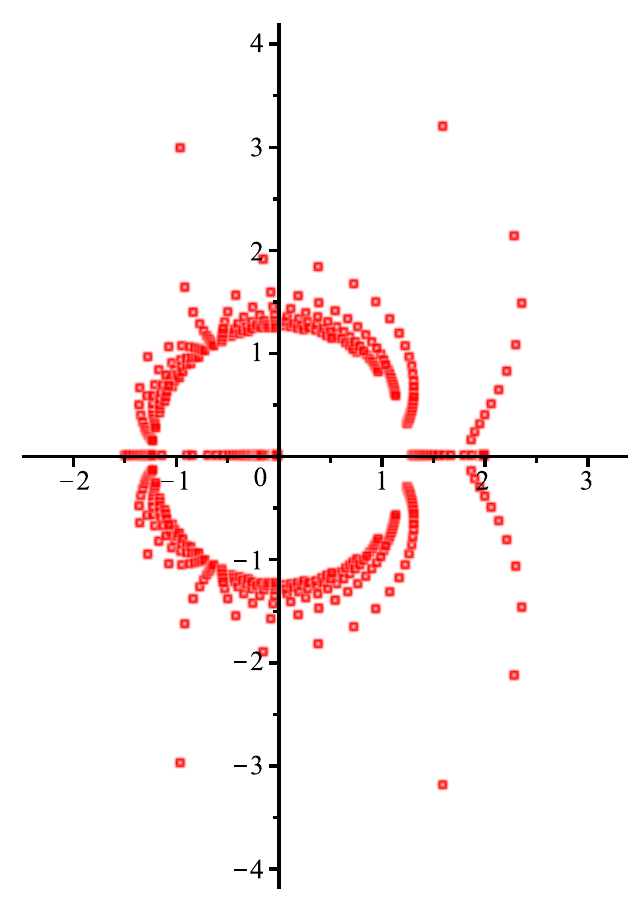} 
  \caption{Zeros of $g_n(x) = x^{n+1} -2 x^n + n^2x^2 +5nx+1$ for $2 \leq n \leq 30$.}
  \end{center}
\end{figure}

We illustrate Theorem~\ref{MainThm} with a couple of examples of just straight sequences of polynomials.  In the next section, we give applications of Theorem~\ref{MainThm} to various graph polynomials. Consider the sequence of polynomials
\[ f_n(x) = x^{n+1} -2 x^n + x^2 +n^2 = (x-2) x^n +(n^2 + x^2)\cdot 1^n\]
By Theorem \ref{MainThm}, the limits of zeros of $f_n$ are the unit circle $|z|=1$ centered at the origin and the isolated point $z = 2$ (see Figure~\ref{example2}). On the other hand,  for
\[ g_n(x) = x^{n+1} -2 x^n + n^2x^2+ 5nx +1 = (x-2) x^n +(n^2x^2+5nx+1)\cdot 1^n,\]
the limits of zeros of $g_n$ are the unit circle centered at the origin and the isolated points $z = 2$ and $z = 0$ (see Figure~\ref{example2}).

\begin{figure}
\label{example1}
\begin{center}
 \includegraphics[scale=0.7]{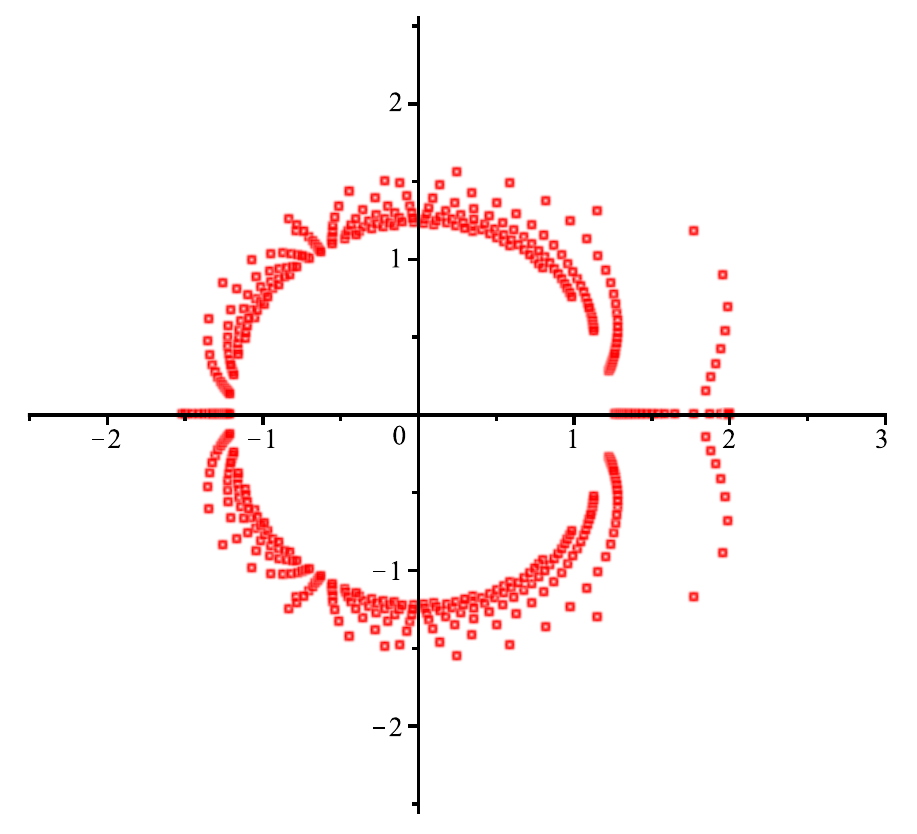} 
 \caption{Zeros of $f_n(x) = x^{n+1} -2 x^n + x^2 +n^2$ for $2 \leq n \leq 30$.}
 \end{center}
\end{figure}

\medskip

\medskip

\section{Applications} 

\medskip

First, to go back to our motivating example of Steele polynomials, Theorem~\ref{MainThm} yields that the limits of zeros of the Steele polynomials for cycle graphs $S(C_n; t)$ defined in \eqref{eqn:SteeleCn} is the unit circle $|z| = 1$.

\begin{figure}
\label{indepexample}
\begin{center}
 \includegraphics[scale=.7]{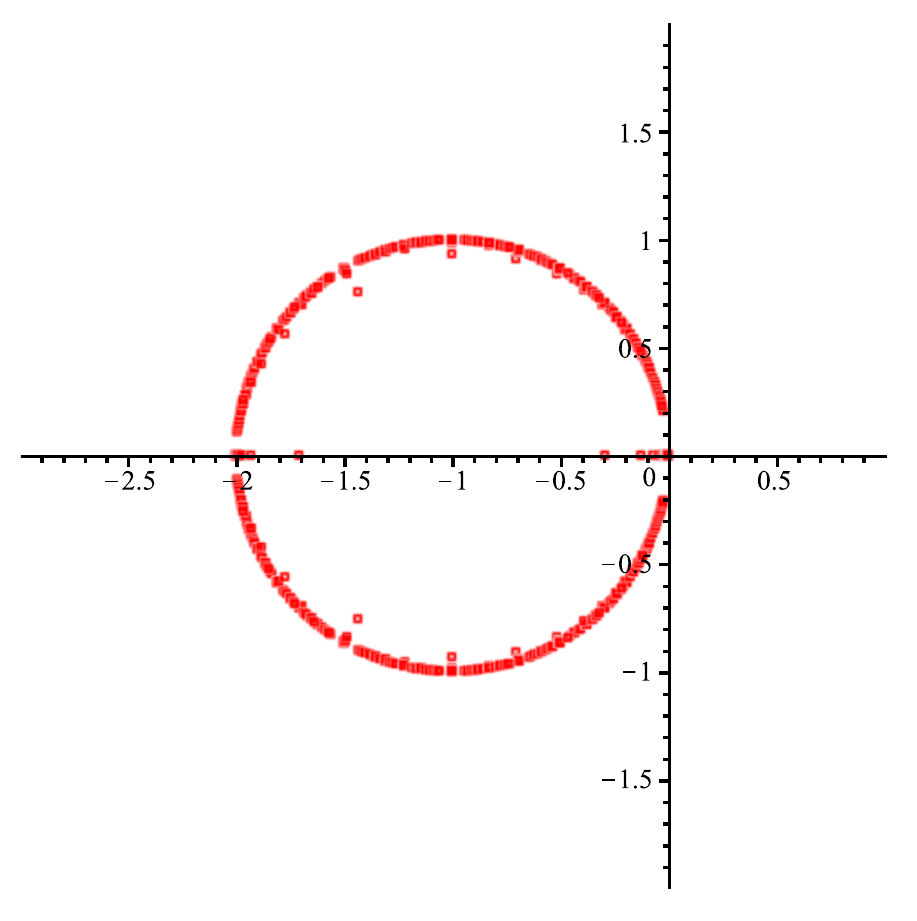} 
 \caption{Zeros of the independence polynomials $n(1+x)^{n} (n-1)$ for $2 \leq n \leq 40$.}
 \end{center}
\end{figure}

We start our list of additional applications with a simple example from independence polynomials. Consider the independence polynomial $i(K_{n,n, \ldots, n},x)$ of a complete $n$-partite graphs $K_{n,n,\ldots,n}$ with $n$ equal parts each of cardinality $n$. As the independent sets are those contained wholly in a part,  $i(K_{n,n,\ldots,n},x) = n(1+x)^n-(n-1)$, which is not of the form for which the original BKW Theorem can be applied. From Theorem~\ref{MainThm} we derive immediately that the limits of zeros are the unit circle $|1+z| = 1$, centered at $z = -1$  (see Figure~\ref{indepexample}). (Of course, the zeros can be found explicitly here as $-1 + \omega_n$, where $\omega_n$ is a an $n$-th zero of $-(n-1)/n$, and these do converge to the circle as $\lim_{n \rightarrow \infty} (n-1)/n = 1$.)

\begin{figure}
 \begin{center}
  \includegraphics[scale=.8]{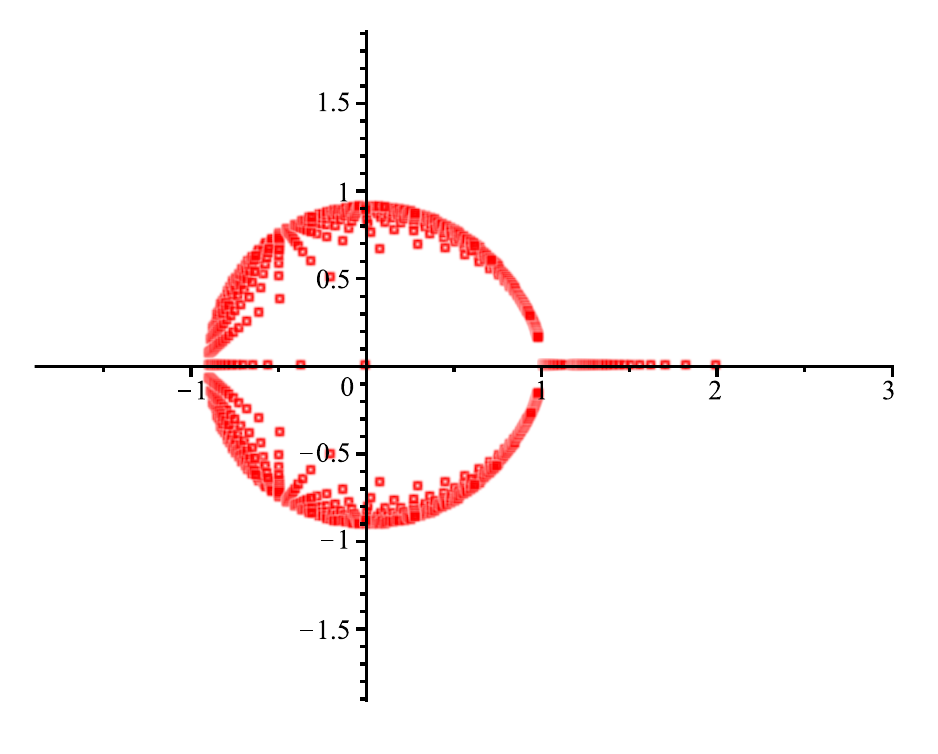} 
 \caption{Zeros of the strongly connected reliability polynomials $2p^{n}-p^{2n}+n(1-p)^2p^{2n-2}$ for $3 \leq n \leq 40$.}
 \end{center}
 \label{screlexample}
\end{figure}

Our second example is afforded from reliability theory. The {\em strongly connected reliability} $\mbox{screl}(D,p)$ (see, for example, \cite{browncox}) is the probability that the spanning subdigraph of directed graph $D$ of operational edges is strongly connected, given that the vertices are always operational but the arcs are independently operational with probability $p$ (this is a natural analogue to all-terminal reliability). Consider the directed graph $C_{n}^{\leftrightarrow}$, formed from the undirected $n$-cycle $C_n$ by replacing each edge $\{k,k+1\}$ by the two arcs $(k,k+1)$ and $(k+1,k)$. It is easy to see that a spanning subdigraph is operational if and only if either one of the two directed cycles is operational, or exactly one pair of parallel arcs is non-operational and all other arcs are operational. It follows (see \cite{browncox}) that, setting $N = n-1$, 
\begin{eqnarray*}
\mbox{screl}(C_{n}^{\leftrightarrow},p) & = & 2p^{n}-p^{2n}+n(1-p)^2p^{2n-2}\\
                                                             & = & 2p^{2} \cdot p^{N} - p^{2} \cdot (p^2)^N + (N+1)(1-p)^2p^{2N}\\
                                                             & = & 2p^{2} \cdot p^{N} + \left( N(1-p)^2 + (1-p)^2 \right) \cdot (p^2)^N
\end{eqnarray*}
Setting $\alpha_1 = 2p^{2}$, $\alpha_2 = N(1-p)^2 + (1-p)^2$, $\lambda_1 = p$, and $\lambda_2 = p^2$, by applying Theorem~\ref{MainThm} we find that the limit of zeros are the unit circle $|p| = 1$ centered at the origin, along with the point $z = 0$ (see Figure~\ref{screlexample}).
           
\begin{figure}
 \begin{center}
 \includegraphics[scale=.8]{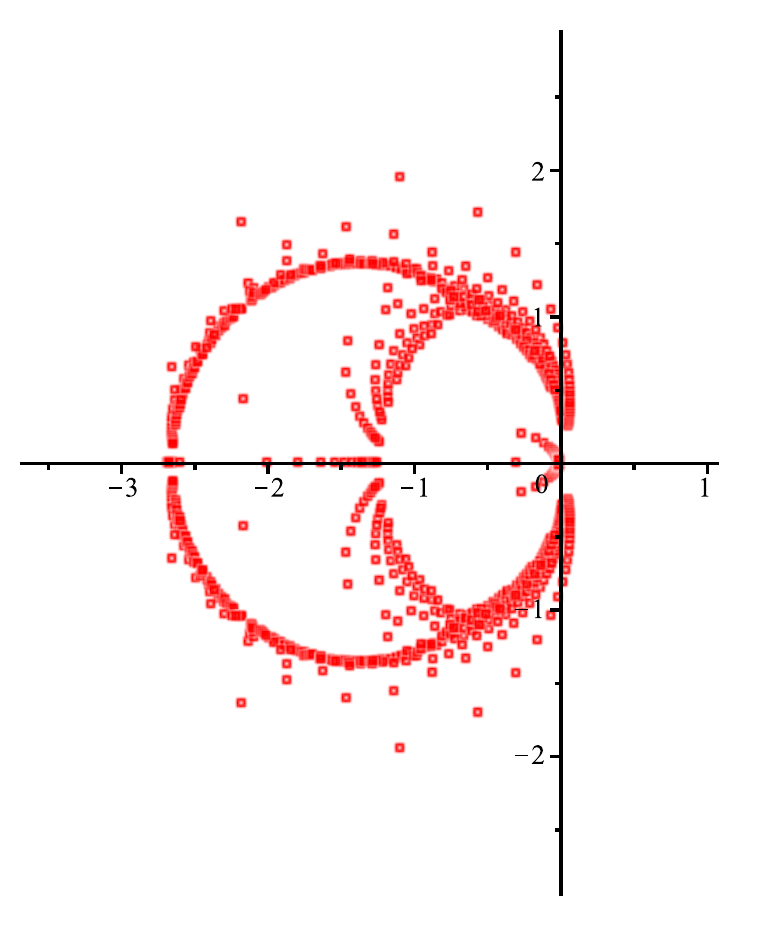} 
 \caption{Zeros of the domination polynomials $((1+x)^n-1-nx)^2 + 2x^n + 2nx^2((1+x)^n-1 - 1)\ + nx^2$ for $2 \leq n \leq 30$.}
  \end{center}
\label{domexample}
\end{figure}

\begin{figure}
 \begin{center}
 \includegraphics[scale=.6]{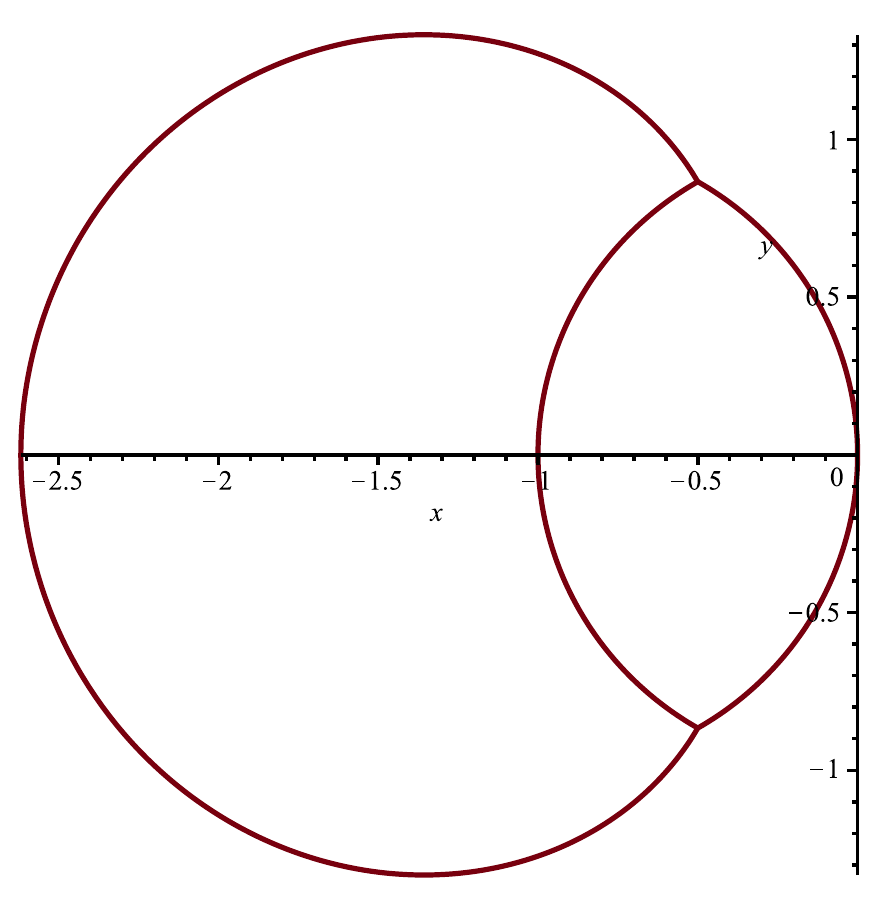} 
 \caption{Limiting curves of the zeros of the domination polynomials $((1+x)^n-1-nx)^2 + 2x^n + 2nx^2((1+x)^n-1 - 1)\ + nx^2$.}
  \end{center}
\label{domcurves}
\end{figure}

Our final example involves domination polynomials For $n \geq 2$, the {\em bipartite cocktail graph} $B_{n}(x)$ is formed from the complete bipartite graph $K_{n,n}$ by removing a perfect matching. The domination polynomial of $B_{n}$ is given by (see \cite{browntufts})  
\[ B_{n}(x) = ((1+x)^n-1-nx)^2 + 2x^n + 2nx^2((1+x)^n-1 - 1)\ + nx^2.\]
We can rewrite this polynomial in the form 
\[ \alpha_1 \lambda_{1}^{N} + \alpha_2 \lambda_{2}^{N} + \alpha_3 \lambda_{3}^{N} + \alpha_4 \lambda_{4}^{N}\]
where
\begin{itemize}
\item $\alpha_1 = (1+x)^2, \lambda_1 = (1+x)^2,$ 
\item $\alpha_2 = (2x^2-2x-2)N + 2x^2-4x-4, \lambda_2 = (1+x),$
\item $\alpha_3 = 2x, \lambda_3  = x,$ and 
\item $\alpha_4 = x^2N+3x^2N+3x^2, \lambda_4 = 1.$
\end{itemize}
It is not hard to check that there are no isolated limits of zeros (as in case (i) of Theorem~\ref{MainThm}). The other limits of zeros come from setting some of the moduli of the $\lambda_i$'s equal and insisting they be larger than the moduli of the rest. Such a case by case study was done in \cite{browntufts} for similar polynomials (all that differs are the $\alpha_i$'s), and the limits of the zeros consist of the union of three curves:    
\begin{itemize}
\item The portion of the unit circle centered at $1$ that has real part at least $-1/2$.
\item The portion of the elliptic curve $|1+z|^2 = |z|$ with real part at most $-1/2$.
\item The portion of the unit circle centered at the origin that has real part at most $-1/2$.
\end{itemize}
Figure~\ref{domexample} shows the zeros of the domination polynomials for $n$ from $2$ to $30$, while Figure~\ref{domcurves} shows the limiting curves.

\section*{Acknowledgements}
 
 
J. Brown acknowledges research support from the Natural Sciences and Engineering Research Council of Canada (NSERC).

\newpage

\bibliographystyle{amsplain}

\begin{thebibliography}{11}

\bibitem{BKW}
S.\ Beraha, J.\ Kahane and N.\ Weiss,
\newblock {\em Limits of zeroes of recursively defined families of polynomials},
\newblock in: Studies in Foundations and Combinatorics, Vol. 1, Advances in Mathematics
Supplementary Studies (G.-C. Rota, ed.), Academic Press, New York
(1978), 213--232.

\bibitem{browncolbourn}
\newblock J.I. Brown and C.J. Colbourn, {\em Roots of the reliability polynomial}, 
\newblock SIAM J. Discrete  Math., {\bf 5} (1992), 571--585.

\bibitem{browncox}
\newblock J.I Brown and D. Cox, {\em The closure of the roots of strongly connected reliability polynomials is the entire complex plane},
\newblock Discrete Math.  {\bf 309} (2009), 5043--5047.

\bibitem{BrownLimit} J.I.\ Brown and C.A.\ Hickman,
\newblock {\em On chromatic roots of large subdivisions of graphs},
\newblock Discrete Mathematics {\bf 242} (2002), 17--30. 

\bibitem{BHN}
\newblock J.I. Brown, C.A. Hickman, and R.J. Nowakowski, {\em On the location of roots  of independence polynomials}, 
\newblock J. Algebraic Combin. {\bf 19} (2004),  273--282.

\bibitem{browntufts} 
\newblock J.I. Brown and J. Tufts, {\em On the Roots of Domination Polynomials}, 
\newblock Graphs  Combin. {\bf 30} (2014), 527--647.

\bibitem{dongbook}
\newblock F.M.\ Dong, , K.M.\ Koh and K.L.\ Teo, {\em Chromatic Polynomials and Chromaticity of Graphs}, World Scientific, London, 2005.

\bibitem{colbook}
C.J. Colbourn, 
\newblock {\em The Combinatorics of Network Reliability}, Oxford University Press, New York, 1987.

\bibitem{KOY} Y.\ Kovchegov, P.\ T.\ Otto, and A.\ Yambartsev,
\newblock {\em Cross multiplicative coalescence and minimal spanning trees of irregular graphs},
\newblock submitted for publication.

\bibitem{NOS} J.\ Nishikawa, P.\ T.\ Otto, and C.\ Starr,
\newblock Polynomial representation  for the expected length of minimal spanning trees,
\newblock Pi Mu Epsilon Journal {\bf 13} (2012), 357--365.

\bibitem{roylesokal}
\newblock G. Royle and A.D. Sokal,
\newblock {\em The Brown-Colbourn conjecture on the zeros of reliability polynomials is false}, 
\newblock J.\ Combin.\ Theory Ser.\ B {\bf 91} (2004),  345--360.

\bibitem{sokaldense}
\newblock A.D. Sokal, {\em Chromatic roots are dense in the whole complex plane}, Probab. Combin. Comput. {\bf 13} (2004), 221--261.

\bibitem{Steele} J.\ M.\ Steele, 
\newblock {\em Minimum spanning trees for graphs with random edge lengths},
\newblock Mathematics and Computer Science {\bf II}  (2002), 223--245.

\end{thebibliography}

\end{document}